\documentclass{article}

\usepackage{indentfirst}
\usepackage{amsmath,amsfonts,amsthm,amssymb}
\usepackage{mathrsfs}

\def\co{\colon\thinspace}
\DeclareMathAlphabet{\mathsfsl}{OT1}{cmss}{m}{sl}

\newtheorem{thm}{Theorem}

\newtheorem{cor}[thm]{Corollary}

\newtheorem*{thm*}{Theorem}

\theoremstyle{definition}
\newtheorem{defn}[thm]{Definition}

\begin{document}

\title{Corrigendum to ``Knot Floer homology detects fibred knots"}

\author{{Yi NI}\\{\normalsize Department of Mathematics, Massachusetts Institute of Technology}\\
{\normalsize 77 Massachusetts Avenue, Cambridge, MA
02139-4307}\\{\small\it Emai\/l\/:\quad\rm yni@math.mit.edu}}

\date{}
\maketitle

\begin{abstract}We correct a mistake on the citation of JSJ theory
in \cite{Ni}. Some arguments in \cite{Ni} are also slightly
modified accordingly.
\end{abstract}

An important step in \cite{Ni} uses JSJ theory \cite{JS,Jo} to
deduce some topological information about the knot complement when
the knot Floer homology is monic, see \cite[Section 6]{Ni}. The
version of JSJ theory cited there is from \cite{CL}. However, as
pointed out by Kronheimer, the definition of ``product regions" in
\cite{CL} is not the one we want. In this note, we will provide
the necessary background on JSJ theory following \cite{JS}. Some
arguments in \cite{Ni} will then be modified.

\begin{defn}
An {\it $n$--manifold pair} is a pair $(M,T)$ where $M$ is an
$n$--manifold and $T$ is an $(n-1)$--manifold contained in
$\partial M$. A $3$--manifold pair $(M,T)$ is {\it irreducible} if
$M$ is irreducible and $T$ is incompressible. An irreducible
$3$--manifold pair $(M,T)$ is {\it Haken} if each component of $M$
contains an incompressible surface.
\end{defn}

%\begin{defn} Let $(M,T)$ be a compact, irreducible $3$--manifold
%pair. A $3$--manifold pair $(\Sigma, \Phi)\subset (M,T)$ is {\it
%perfectly-embedded} in $(M,T)$ if (i) $\Sigma\cap \partial
%M=\Phi$, (ii) $\mathrm{Fr(\Sigma)}$ is incompressible in $M$,
%(iii) no component of $\mathrm{Fr(\Sigma)}$ is $T$--parallel, (iv)
%$(\Sigma,\Phi)$ has no
%\end{defn}

\begin{defn}\cite[Page 10]{JS}
A compact $3$--manifold pair $(\mathcal S,\mathcal T)$ is called
an {\it$I$--pair} if $\mathcal S$ is an $I$--bundle over a compact
surface, and $\mathcal T$ is the corresponding $\partial
I$--bundle. A compact $3$--manifold pair $(\mathcal S,\mathcal T)$
is called an {\it$S^1$--pair} if $\mathcal S$ is a Seifert fibred
manifold and $\mathcal T$ is a union of Seifert fibres in some
Seifert fibration. A {\it Seifert pair} is a compact $3$--manifold
pair $(\mathcal S,\mathcal T)$, each component of which is an
$I$--pair or an $S^1$--pair.
\end{defn}

\begin{defn}\cite[Page 138]{JS}
A {\it characteristic pair} for a compact, irreducible
$3$--manifold pair $(M,T)$ is a perfectly-embedded Seifert pair
$(\Sigma,\Phi)\subset (M,\mathrm{int}(T))$ such that if $f$ is any
essential, nondegenerate map of an arbitrary Seifert pair
$(\mathcal S,\mathcal T)$ into $(M,T)$, $f$ is homotopic, as a map
of pairs, to a map $f'$ such that $f'(\mathcal S)\subset\Sigma$
and $f'(\mathcal T)\subset\Phi$.
\end{defn}

The definition of a perfectly-embedded pair can be found in
\cite[Page 4]{JS}. We note that the definition requires that
$\Sigma\cap\partial M=\Phi$, so $\Sigma$ is disjoint from
$\partial M-T$.

The main result in JSJ theory is the following theorem.

\begin{thm}[Jaco--Shalen \cite{JS}, Johannson \cite{Jo}]
Every Haken $3$--manifold pair $(M,T)$ has a characteristic pair.
This characteristic pair is unique up to ambient isotopy relative
to $(\partial M-\mathrm{int}(T))$.
\end{thm}

%For our purpose, we will introduce some other notions.

\begin{defn}
Let $(M,\gamma)$ be a sutured manifold. A $3$--manifold pair
$(P,Q)\subset (M,R(\gamma))$ is a {\it product pair} if
$P=F\times[0,1], Q=F\times\{0,1\}$ for some compact surface $F$,
and $F\times0\subset R_-(\gamma), F\times1\subset R_+(\gamma)$. We
also require that $P\cap A=\emptyset \;\text{or}\;A$ for any
annular component $A$ of $\gamma$. A product pair is {\it gapless}
if no component of its exterior is a product pair.
\end{defn}

\begin{defn}
Suppose $(M,\gamma)$ is a taut sutured manifold, $(\Sigma,\Phi)$
is the characteristic pair for $(M,R(\gamma))$. The {\it
characteristic product pair} for $M$ is the union of all
components of $(\Sigma,\Phi)$ which are product pairs. A {\it
maximal product pair} for $M$ is a gapless product pair $(\mathcal
P,\mathcal Q)$ such that it contains the characteristic product
pair, and if $(\mathcal P',\mathcal Q')\supset(\mathcal P,\mathcal
Q)$ is another gapless product pair, then there is an ambient
isotopy relative to $\gamma$ that takes $(\mathcal P',\mathcal
Q')$ to $(\mathcal P,\mathcal Q)$.
\end{defn}

The existence of maximal product pairs follows from the
definition, although the uniqueness is not guaranteed. The
exterior of a maximal product pair is also a sutured manifold. By
definition the exterior does not contain essential product annuli
or essential product disks.

Now we are ready to modify the arguments in \cite{Ni}. The next
theorem is a reformulation of \cite[Theorem 6.2]{Ni}. The proof is
not changed though.

\vspace{5pt}

\noindent{\bf Theorem 6.2$'$} {\it Suppose $(M,\gamma)$ is an
irreducible balanced sutured manifold, $\gamma$ has only one
component, and $(M,\gamma)$ is vertically prime. Let $\mathcal E$
be the subgroup of $H_1(M)$ spanned by the first homologies of
product annuli in $M$. If $\widehat{HFS}(M,\gamma)\cong\mathbb Z$,
then $\mathcal E=H_1(M)$.}\qed

\vspace{5pt}

\begin{cor}\label{CharSurj}
In the last theorem, suppose $(\Pi,\Psi)$ is the characteristic
product pair for $M$, then the map
$$i_*\co H_1(\Pi)\to H_1(M)$$
is surjective.
\end{cor}
\begin{proof}
We recall that such an $M$ is a homology product
\cite[Proposition~3.1]{Ni}.

Suppose $(\Sigma,\Phi)$ is the characteristic pair for
$(M,R(\gamma))$, then any product annulus can be homotoped into
$(\Sigma,\Phi)$ without crossing $\gamma$. Let
$\Phi_+=\big(\Phi\cap
R_+(\gamma)\big)\subset\mathrm{int}(R_+(\gamma))$. Theorem~6.2$'$
implies that the map $H_1(\Phi_+)\to H_1(R_+(\gamma))$ is
surjective, so $\partial\Phi_+$ consists of separating circles in
$R_+(\gamma)$. If a component $(\mathcal S,\mathcal T)$ of
$(\Sigma,\Phi)$ is an $S^1$--pair, then $\mathcal T\cap
R_+(\gamma)$ consists of annuli by definition. We conclude that
each annulus is null-homologous in $H_1(R_+(\gamma))$.

Suppose a product annulus $A$ contributes to $H_1(M)$
nontrivially, and it can be homotoped into a component
$(\sigma,\varphi)$ of $(\Sigma,\Phi)$. Given the result from the
last paragraph, this $(\sigma,\varphi)$ cannot be an $S^1$--pair.
It is neither a twisted $I$--bundle since the two components of
$\partial A$ are contained in different components of $R(\gamma)$.
So $(\sigma,\varphi)$ must be a trivial $I$--bundle, and the two
components of $\varphi$ lie in different components of
$R(\gamma)$. In other words, $(\sigma,\varphi)$ is a product pair.
Now our desired result follows from Theorem~6.2$'$.
\end{proof}

The following proof of the main theorem in \cite{Ni} is only
slightly changed. Basically we use ``maximal product pair" here
instead of the wrong notion ``characteristic product region" in
\cite{Ni}.

\vspace{5pt}

\noindent{\it Proof of {\rm\cite[Theorem 1.1]{Ni}}.}\quad Suppose
$(M,\gamma)$ is the sutured manifold obtained by cutting open
$Y-\mathrm{int}(\mathrm{Nd}(K))$ along $F$, $(\mathcal P,\mathcal
Q)$ is a maximal product pair for $M$. We need to show that $M$ is
a product. By \cite[Proposition~3.1]{Ni}, $M$ is a homology
product. Moreover, by \cite[Theorem~4.1]{Ni}, we can assume $M$ is
vertically prime.

If $M$ is not a product, then $M-\mathcal P$ is nonempty. Thus
there exist some product annuli in $(M,\gamma)$, which split off
$\mathcal P$ from $M$. Let $(M',\gamma')$ be the remaining sutured
manifold. By definition $(\mathcal P,\mathcal Q)$ contains the
characteristic product pair for $M$. Corollary~\ref{CharSurj} then
implies that the map $H_1(\mathcal P)\to H_1(M)$ is surjective. So
$R_{\pm}(\gamma')$ are planar surfaces, and $M'\cap\mathcal P$
consists of separating product annuli in $M$. Since we assume that
$M$ is vertically prime, $M'$ must be connected. (See the first
paragraph in the proof of \cite[Theorem~5.1]{Ni}.) Moreover, $M'$
is also vertically prime. By \cite[Theorem~5.1]{Ni},
$\widehat{HFS}(M',\gamma')\cong\mathbb Z$.

We add some product 1--handles to $M'$ to get a new sutured
manifold $(M'',\gamma'')$ with $\gamma''$ connected. By
\cite[Proposition~2.9]{Ni},
$\widehat{HFS}(M'',\gamma'')\cong\mathbb Z$. It is easy to see
that $M''$ is also vertically prime. \cite[Proposition~3.1]{Ni}
shows that $M''$ is a homology product.

Let $H$ be one of the product $1$--handles added to $M'$ such that
$H$ connects two different components of $\gamma'$. By
Theorem~6.2$'$, there is at least one product annulus $A$ in
$M''$, such that $A$ cannot be homotoped to be disjoint from the
cocore of $H$. Isotope $A$ if necessary, we find that at least one
component of $A\cap M'$ is an essential product disk in $M'$, a
contradiction to the assumption that $(\mathcal P,\mathcal Q)$ is
a maximal product pair. \qed

\vspace{5pt}

\noindent{\bf Acknowledgements.}\quad We are extremely grateful to
Peter Kronheimer for pointing out the mistake in \cite{Ni}. This
note is written when the author visited Zhejiang University. The
author wishes to thank Feng Luo for his hospitality during the
visit.

\end{document}